\documentclass[10pt]{amsart}

\usepackage{amssymb,latexsym,amsmath,amsfonts}
\usepackage[mathscr]{eucal}
\usepackage{mathrsfs}
\usepackage{pxfonts}
\usepackage{amscd}
\usepackage{amsthm}
\usepackage{amsxtra}
\usepackage[nointegrals]{wasysym}
\usepackage{bm}
\usepackage{calc}
\usepackage{float}
\usepackage[usenames,dvipsnames]{xcolor}
\usepackage{hyperref}
\hypersetup{colorlinks=true, citecolor=MidnightBlue, linkcolor=BrickRed}
 \usepackage[normalem]{ulem}

\numberwithin{equation}{section}
\theoremstyle{definition}
\newtheorem{definition}{Definition}[section]
\theoremstyle{definition}

\theoremstyle{plain}
\newtheorem{theorem}[definition]{Theorem}

\newcommand{\beas}{\begin{eqnarray*}}
\newcommand{\eeas}{\end{eqnarray*}}
\newcommand{\bes} {\begin{equation*}}
\newcommand{\ees} {\end{equation*}}
\newcommand{\be} {\begin{equation}}
\newcommand{\ee} {\end{equation}}
\newcommand{\bea} {\begin{eqnarray}}
\newcommand{\eea} {\end{eqnarray}}

\newcommand{\zt}{\zeta}

\newcommand{\om}{\omega}

\newcommand{\dbar}{\overline\partial}

\newcommand{\Om}{\Omega}

\newcommand{\bh}{\boldsymbol h}

\newcommand{\zbar}{\overline z}

\newcommand{\wt}{\widetilde}

\newcommand{\hol}{\mathcal{O}}

\newcommand{\C} {\mathbb{C}}

\newcommand{\p}{\mathbb{P}^1}

\begin{document}

\title[Dimension of Bergman spaces]{On the dimension of bundle-valued Bergman spaces on compact Riemann surfaces}

\author{A.-K. Gallagher}
\address{Gallagher Tool \& Instrument LLC, Redmond, WA}
\email{anne.g@gallagherti.com}
\author{P. Gupta}
\address{Department of Mathematics, Indian Institute of Science, Bangalore}
\email{purvigupta@iisc.ac.in}
\author{L. Vivas}
\address{Department of Mathematics, The Ohio State University, Columbus, OH}
\email{vivas@math.osu.edu}

\dedicatory{Dedicated to the memory of Berit Stens{\o}nes}


\begin{abstract}
Given a holomorphic vector bundle $E$ over a compact Riemann surface $M$, and an open set $D$ in $M$, we prove that the Bergman space of holomorphic sections of the restriction of $E$ to $D$ must either coincide with the space of global holomorphic sections of $E$, or be infinite dimensional. Moreover, we characterize the latter entirely in terms of potential-theoretic properties of $D$. 
\end{abstract}

\maketitle

\section{Introduction}\label{S:Intro}


Let $(M,g)$ be a compact Riemann surface equipped with a Hermitian metric, and $(E,h)$ be a Hermitian holomorphic vector bundle over $M$. Given an open subset $D\subset M$, the Bergman space of $E$-valued holomorphic sections on $D$ is the reproducing kernel Hilbert space given by
	\bes
		A^2(D,g;E,h)=\left\{s\in\hol(D;E):||s||_h:=\left(\int_D h(s,s)\,\omega_g\right)^{1/2}<\infty\right\},
	\ees
where $\hol(D;E)$ is the space of holomorphic sections of $E$ over $D$, and $\om_g$ is the volume form associated to $g$. As a vector space, $A^2(D;E)=A^2(D,g;E,h)$ is independent of the choice of $g$ and $h$. We prove the following theorem. 

\begin{theorem}\label{T:main} Let $M$ be a compact Riemann surface and $E$ be a holomorphic vector bundle over $M$. Suppose $D\subset M$ is a nonempty open subset of $M$. Then the following are equivalent.
	\begin{itemize}
\item[$(a)$] $M\setminus D$ is nonpolar.
\item[$(b)$] $\hol(M;E)\subsetneq A^2(D;E)$.
\item[$(c)$] $\dim A^2(D;E)=\infty$. 
\item[$(d)$] There exists a bounded subharmonic function $\psi\in\mathcal{C}^\infty(D)$ such that $i\partial\dbar\psi\geq \omega$ on $D$ for some volume form $\omega$ on $M$.
\end{itemize}  
\end{theorem}  

This settles the question posed by Sz{\H o}ke in \cite{Sz20} on the dimension of vector bundle-valued Bergman spaces on open subsets of a compact Riemann surface. Sz{\H o}ke's question is motivated by the well-understood case of open subsets of the complex plane. A complete analogue of Theorem ~\ref{T:main} for open subsets of $\C$ holds by the works of Carleson \cite{Ca67}, Wiegerinck \cite{Wi84}, Gallagher-Harz-Herbort \cite{GaHaHe17}, and Gallagher-Lebl-Ramachandran \cite{GaLeRa21}. Using an argument similar to Wiegerinck's, Sz{\H o}ke established the equivalence of $(b)$ and $(c)$ in Theorem \ref{T:main} for the special case of $M=\p$. A proof of the implication $(b)\Rightarrow (a)$ for a general $M$ is also given in \cite{Sz20}. The rest of the theorem for $M=\p$ is established by the authors in \cite{GaGuVi22}, using the analytic approach of \cite{GaHaHe17} and \cite{GaLeRa21}. In comparison to Sz{\H o}ke's argument, the latter approach is not as sensitive to the topological and algebraic structure of the underlying surface and bundle. Thus, it lends itself more easily to generalization to all compact Riemann surfaces. 

Since $(c)$ trivially implies $(b)$, and it is proven in \cite{Sz20} that $(b)$ implies $(a)$, the proof of Theorem~\ref{T:main} reduces to proving that $(a)$ implies $(d)$, and $(d)$ implies $(c)$. For the former, the Green's function of the domain is suitably modified to obtain the desired strictly subharmonic function; see Section~\ref{SS:(a)=>(d)}. For the latter, one relies on a version of H{\"o}rmander's $L^2$-method for solving the Cauchy--Riemann equations for bundle-valued forms; see Theorem~\ref{T:HDtheorem}. This cannot be directly invoked for the given Hermitian metric $h$ on $E$, as its curvature may not be positive. We instead apply H{\"o}rmander's theorem to the twisted metric $e^{-\wt\psi}h$, where $\wt\psi$ is a nonsmooth correction of $\psi$ as in $(d)$, chosen so as to yield linearly independent holomorphic sections in the untwisted Bergman space. Note that the twisted metric is singular, but the $L^2$-methods cited here are for smooth metrics. Due to the benign nature of the singularities, this issue is resolved by using an extension trick that appears to be standard in this context; see \cite[Lemma 5.1.3]{Be10}, for instance.

\medskip
\noindent{\bf Acknowledgements.} This material is partially based upon work supported by the National Science Foundation under Grant No. DMS-1928930 while the authors participated in the Summer Resarch in Mathematics program hosted
by the Mathematical Sciences Research Institute in Berkeley, California, in June 2022. The authors are grateful to MSRI for providing financial support and a stimulating environment. Gupta and Vivas were partially supported by an Infosys Young Investigator Award, and NSF Grant No. DMS-1800777, respectively.

\section{Background and preliminaries}\label{S:prelim} 
\subsection{Holomorphic vector bundles on compact Riemann surfaces} Throughout this section, $X$ is a Riemann surface equipped with a Hermitian metric $g$. The holomorphic and antiholomorphic tangent bundles of $X$ are denoted by $T^{1,0}X$ and $T^{0,1}X$, respectively. For $p,q\in\{0,1\}$, ${T^*}^{p,q}X$ is the bundle of $(p,q)$-forms on $X$ and $\Om^{p,q}(X)$ is the space of its smooth sections. The standard Dolbeault operators on $\Om^{p,q}(X)$ are denoted by $\partial$ and $\dbar$. 

 The volume form associated to $g$ is denoted by $\om_g$. Note that $g$ and $\overline g$ are Hermitian metrics on $T^{1,0}X$ and $T^{0,1}X$, respectively. The dual Hermitian metrics on ${T^*}^{1,0}X$ and ${T^*}^{0,1}X$ are denoted by $g^*_{1,0}$ and $g^*_{0,1}$, respectively. The complex bundle ${T^*}^{1,1}X$ is endowed with the metric $g^*_{1,1}$ given on ${T^*}^{1,1}_xX$ by 
	\bes
		g^*_{1,1}\left(\sigma\wedge\tau,\sigma'\wedge\tau' \right)=g^*_{1,0}(\sigma,\sigma')g^*_{0,1}(\tau,\tau'),
	\ees
for $\sigma,\sigma'\in {T^*}^{1,0}_xX$ and $\tau,\tau'\in {T^*}^{0,1}_xX$, and  extended to general multilinear forms by sesquilinearity. Lastly, $g^*_{0,0}$ is the metric on ${T^*}^{0,0}X\cong X\times\C$ given at each $x\in X$ by the standard Hermitian metric on $\C$. In local coordinates, if $g=\gamma\: dz\otimes d\zbar$, then $\om_g=i\:\gamma\: dz\wedge d\zbar$, $g^*_{1,0}(adz,bdz)=\frac{a\overline b}{\gamma}$ and $g^*_{1,1}(adz\wedge d\zbar,bdz\wedge d\zbar)=\frac{a\overline b}{\gamma^2}$.

Given a holomorphic vector bundle $E$ over $X$, its spaces of smooth and holomorphic sections are denoted by $\Gamma(X;E)$ and $\hol(X;E)$, respectively. An $E$-valued $(p,q)$-form is a section of ${T^*}^{p,q}X\otimes E$, and for convenience, $\Gamma(X;{T^*}^{p,q}X\otimes E)$ is denoted by $\Om^{p,q}(X;E)$. In particular, $\Om^{0,0}(X;E)=\Gamma(X;E)$. The space of smooth compactly supported $E$-valued $(p,q)$-forms is denoted by $\Om_c^{p,q}(X;E)$. For $p\in\{0,1\}$, the operator $\dbar:\Om^{p,0}(X;E)\rightarrow \Om^{p,1}(X;E)$ is given in local frames by 
	\bes
		\dbar:s=\sum_{j=1}^r \sigma_j\otimes \xi_j\mapsto \sum_{j=1}^r \dbar\sigma_j\otimes \xi_j,
	\ees
where $\sigma_1,...,\sigma_r\in \Om^{p,0}(X)$ and $\xi_1,...,\xi_r$ is a choice of local frames of $E$.
Given a Hermitian metric $h$ on $E$, ${T^*}^{p,q}X\otimes E$ is endowed with the Hermitian metric ${g^*}_{\!\!p,q}\otimes h$. This induces the following inner product structure on $\Om^{p,q}(X;E)$:
	\beas
		\left<s,t\right>_{g^*_{p,q}\otimes h}
		&=&\int_X(g^*_{p,q}\otimes h)(s,t)\:\om_g,\qquad s,t\in\Om^{p,q}(X;E).
	\eeas
For simplicity, we denote $\left<s,t\right>_{g^*_{p,q}\otimes h}$ by $\left<s,t\right>_{h}$ since $g$ is fixed throughout. The Hilbert space completion of the inner product space $(\Om_c^{p,q}(X;E),\left<\cdot,\cdot\right>_h)$ is denoted by $L^2_{p,q}(X,g;E,h)$. One extends $\dbar$ to an operator on $L^2_{p,0}(X,g;E,h)$, $p\in\{0,1\}$, in the sense of currents via the maximal extension of $\dbar$; see \cite[Chapter VIII, \S 3]{De97}. Recall that the Bergman space of sections of $(E,h)$ is defined as $A^2(X,g;E,h)=\hol(X;E)\cap L^2_{0,0}(X,g;E,h)$. If $X$ is compact and $D\subset X$ is an open subset, the vector space $A^2(D,g;E,h)$ is independent of the choice of Hermitian metrics on $X$ and $E$; see Lemma 2.2 in \cite{GaGuVi22}, for instance. Thus, $A^2(D,g;E,h)$ and  $A^2(D,g;{T^*}^{1,0}X\otimes E,g^*_{1,0}\otimes h)$ are simply denoted by $A^2(D;E)$ and $A^2_{1,0}(D;E)$, respectively.

For the rest of this section, $X$ is assumed to be noncompact. Since every holomorphic vector bundle over a noncompact Riemann surface is trivial, see \cite[\S 8.2]{Varolin11}, we may fix a global holomorphic frame, say $\xi=\{\xi_1,...,\xi_r\}$, of $E$ for computational ease. Any $s\in \Om^{p,q}(X;E)$ can then be written uniquely as 
	\bes
		s=\sum_{j=1}^r \sigma_j\otimes \xi_j,
	\ees
for some $\sigma_1,...,\sigma_r\in\Om^{p,q}(X)$. The Hermitian metric $h$ is represented by the positive-definite Hermitian matrix-valued smooth function $h=\left(h_{jk}\right)$ on $X$, where $h_{jk}=h(\xi_j,\xi_k)$, $j,k\in\{1,...,r\}$. The Chern curvature tensor of $(E,h)$ is the $\text{Herm}(E,E)$-valued $(1,1)$-form  $\Theta(h)=\dbar\left(\overline h^{-1}\partial \overline h\right)$. It may be expressed as the matrix of $(1,1)$-forms given by
	\bes
		\Theta(h)=\dbar\left(h^{\ell j}\partial h_{k\ell} \right),
	\ees
where $(h^{jk})$ is the inverse of $h=(h_{jk})$. By further choosing a global holomorphic frame $Z$ of ${T^*}^{1,0}X$, we may also write $\Theta(h)=\left(\theta_{jk}\right) Z\wedge\overline Z$, for some Hermitian-matrix valued smooth function $\theta=\left(\theta_{jk}\right)$ on $X$. Given two $\text{Herm}(E,E)$-valued $(1,1)$-forms $A=\left(a_{jk}\right) Z\wedge\overline Z$ and $B=\left(b_{jk}\right) Z\wedge\overline Z$, we say that $A\geq B$, if $(a_{jk})-(b_{jk})$ is positive semi-definite. This notion is independent of the choice of the global holomorphic frame $Z$. 

Finally, in order to produce holomorphic functions on a noncompact Riemann surface, we will use the following version of H{\"o}rmander's theorem for the $\dbar$-problem on bundle-valued forms. 

\begin{theorem}\cite[Theorem 6.1]{De97}\label{T:HDtheorem} Let $X$ be a noncompact Riemann surface with Hermitian metric $g$ and volume form $\om_g$. Let $(F,\bh)$ be a Hermitian holomorphic vector bundle over $X$ such that
\be\label{E:positive}
i\Theta(\bh)\geq c\om_g\otimes\text{I}_F
\ee
for some $c>0$. Then, for every $\alpha\in L^2_{1,1}(X,g;F,\bh)$,
there exists a $u\in L^2_{1,0}(X,g;F,\bh)$ such that $\dbar u=\alpha$ in the sense of distributions, and
	\bes
		\left<u,u\right>_{\bh}\leq \frac{1}{c}\left<\alpha,\alpha\right>_{\bh}.
	\ees
\end{theorem}

Note that Theorem 6.1 in \cite{De97} is stated more generally for an $m$-semi-positive vector bundle $(F,\bh)$ over a complete K{\"a}hler $n$-dimensional manifold $(X,\hat g)$ equipped with a possibly noncomplete metric $g$. Moreover, $\Om_c^{1,1}(X;F)$ is endowed with the inner product: $(\alpha,\beta)=\int_X g^*_{1,1}(A^{-1}\alpha,\beta)\om_g$, where $A=i\Theta(\bh)\wedge\Lambda$ and $\Lambda$ is given by $g^*_{0,0}(\Lambda\alpha,f)=g^*_{1,1}(\alpha,f\wedge\om_g)$, $\alpha\in\Om^{1,1}(X;F)$, $f\in\Om^{0,0}(X;F)$. In our case, since $X$ is a noncompact Riemann surface, it admits a complete K{\"a}hler metric; see \cite[Theorem~5.2]{De97}. Furthermore, since $m=n=1$, \eqref{E:positive} is equivalent to both the $m$-positivity of $(F,\bh)$ and the positive-definiteness of $A$. For the latter, observe that when $n=1$, the operator $A:\Om^{1,1}(X;F)\rightarrow \Om^{1,1}(X;F)$ acts as $\sigma\otimes\xi\mapsto \frac{1}{\gamma}\sigma\otimes \Theta(\bh)(\xi)$, $\sigma\in\Om^{1,1}(X)$, $\xi\in\Gamma(X;F)$, where we view $\Theta(\bh)$ as an $\text{End}(F)$-valued $(1,1)$-form on $X$.


\subsection{Potential theory on Riemann surfaces} Let $X$ be a Riemann surface. An upper semicontinuous function $s:X\rightarrow [-\infty, \infty)$ with $s\nequiv -\infty$ on any connected component of $X$, is said to be {\em subharmonic on $X$} if, for every coordinate chart $(V,\varphi)$ on $X$, $s\circ\varphi^{-1}$ is subharmonic on $\varphi(V)\subset \C$. The notion of harmonicity is defined analogously via charts. Of particular importance are the so-called Green's functions on $X$. \begin{definition}\label{D:Greensfunction}
  Given $x\in X$, a {\em Green's function on $X$ with singularity at $x$} is a subharmonic function $G_x:X\rightarrow[-\infty,0)$ such that 
	\begin{itemize}
    \item [$(i)$] $G_x$ is harmonic on $X\setminus\{x\}$,
    \item [$(ii)$] if $(U,\varphi)$ is a coordinate chart containing $x$, then $y\mapsto G_x(y)-\log|\varphi(x)-\varphi(y)|$ is harmonic on $U$,  
    \item [$(iii)$] if $H:X\rightarrow[-\infty,0)$ is a subharmonic function satisfying $(i)$ and $(ii)$, then $H\leq G_x$.
  \end{itemize}  
  If $G_x$ exists for every $x\in X$, we call $G(x,y):=G_x(y)$, $x,y\in X$, the {\em Green's kernel} of $X$, and say that $X$ admits a Green's kernel.
\end{definition}

For an open subset of a compact Riemann surface, the existence of the Green's kernel is completely characterized by the nonpolarity of its complement. Recall that a set $K\subset X$ is said to be {\em polar} if for each $\zt\in X$, there is an open neighborhood $V_\zt\subset X$ of $\zt$ and a subharmonic function $s_\zt$ on $V_\zt$ such that $K\cap V_\zt\subset\{x\in V_\zt:s_\zt(x)=-\infty\}$. We recall the following well-known result.

\begin{theorem}\label{T:Greenian} Let $X$ be a compact Riemann surface and $D\subset X$ be an open subset of $X$. Then the following are equivalent.
	\begin{itemize}
\item [$(a)$] $X\setminus D$ is nonpolar.
\item [$(b)$] $D$ admits a Green's kernel. 
\item [$(c)$] $D$ is potential-theoretically hyperbolic, i.e., $D$ admits a nonconstant bounded subharmonic function. 
\end{itemize} 
\end{theorem}
Since the above theorem is a combination of several known results, we give some references. That $(a)$ implies $(c)$ follows from Kellog's theorem,  see \cite[Theorem 4.2.5]{Ra95}, and \cite[Lemma 7.1.5]{Varolin11}. For $(c)\Rightarrow(b)$; see Theorem 7.1.13 in \cite{Varolin11}. Finally, the removability of closed polar sets for bounded subharmonic functions gives that $(b)$ implies $(a)$.

\section{Proof of Theorem~\ref{T:main}}\label{S:proofs}
Since $\hol(M;E)$ is finite dimensional whenever $M$ is compact, $(c)\Rightarrow(b)$. For the proof of $(b)\Rightarrow(a)$; see \cite[Proposition 2.1]{Sz20}. In Sections~\ref{SS:(a)=>(d)} and \ref{SS:dbar}, we prove the implications $(a)\Rightarrow(d)$ and $(d)\Rightarrow(c)$, respectively. 

\subsection{Proof of $\mathbf{(a)}$ implies $\mathbf{(d)}$}\label{SS:(a)=>(d)}

Suppose $M\setminus D$ is nonpolar. Then, by Theorem \ref{T:Greenian}, $D$ admits a Green's kernel $G$. Fix an $x_0\in D$, 
and let $G_{x_0}(y):=G(x_0,y)$ be the associated Green's function on $D$. Define $\psi_{0}:D\longrightarrow [0,1)$ by
\[ 
  \psi_0(y):=
  \begin{cases} 
      e^{2G_{x_0}(y)}, &\text{if }  y\neq x_0, \\
      0, &\text{if } y=x_0.
   \end{cases}
\]
We first show that $\psi_0\in\mathcal{C}^\infty(D)$ is a subharmonic function which is strictly subharmonic outside a discrete set $Z$. By definition of the Green's function, see $(i)$ of Definition \ref{D:Greensfunction}, smoothness of $\psi_0$ only needs to be verified at $y=x_0$. For that, let $(U,\varphi)$ be a coordinate chart of $D$ which includes $x_0$. Then, by $(ii)$ of Definition \ref{D:Greensfunction} there exists a harmonic function $h$ on $U$ such that
$$G_{x_0}(y)-\ln|\varphi(x_0)-\varphi(y)|=h(y)\qquad\forall y\in U.$$
Therefore,
\begin{align}\label{E:smoothpsi_0}
  \psi_0(y)=e^{2h(y)}|\varphi(x_0)-\varphi(y)|^2\qquad\forall y\in U.
\end{align}
It follows that $\psi_0\in\mathcal{C}^\infty(D)$. To show that $\psi_0$ is subharmonic on $D$, we compute
for $y\in D\setminus\{x_0\}$ that
\begin{align}\label{E:subharmonicpsi_0}
  i\partial\overline{\partial}\psi_0(y)
  &=e^{2G_{x_0}(y)}\left(4i\partial G_{x_0}(y)\wedge\overline{\partial}G_{x_0}(y)+2i\partial\overline{\partial}G_{x_0}(y)\right)\notag\\
  &=e^{2G_{x_0}(y)}\left(4i\partial G_{x_0}(y)\wedge\overline{\partial}G_{x_0}(y)\right),
\end{align}
where the last step follows from the harmonicity of $G_{x_0}$ on $D\setminus\{x_0\}$. A similar computation, using  identity \eqref{E:smoothpsi_0} and the holomorphicity of $\varphi$, yields
\begin{align*}
  i\partial\overline{\partial}\psi_0(x_0)=\left(e^{2h}i\partial\varphi\wedge\overline{\partial\varphi}\right)(x_0).
\end{align*}
Thus, $\psi_0$ is strictly subharmonic near $x_0$. It then follows from \eqref{E:subharmonicpsi_0}, that $\psi_0$ is strictly subharmonic on $D\setminus Z$ with
$$Z=\{y\in D\setminus\{x_0\}:\partial G_{x_0}(y)=0 \}.$$
Note that the harmonicity of $G_{x_0}$ implies that $\partial G_{x_0}$ is a holomorphic $(1,0)$-form. Hence, $Z$ is a discrete set in $D\setminus\{x_0\}$. To accommodate for the lack of strict subharmonicity on $Z$ and uniform strict subharmonicity near $bD$, we construct a particular set of cut-off functions. 

We may assume that $M\setminus D$ lies in a single coordinate chart. If not, let $(U,\varphi)$ be a coordinate chart such that $U\setminus D$ is nonpolar. Let $K\subset U\setminus D$ be a nonpolar compact subset. Then $\wt D:=M\setminus K$ is an open set in $M$ containing $D$. Thus, if $(d)$ of Theorem \ref{T:main} holds on $\wt D$, it holds on $D$.

Let $(U_1,\varphi_1)$ be the coordinate chart which contains $M\setminus D$, and $V_1\Subset U_1$ be an open neighborhood of $M\setminus D$. 
Then, $Z\setminus V_1$ is a finite set of points $\{x_2,\dots,x_m\}$. 
For each $j\in\{2,\dots,m\}$, choose a coordinate chart $(U_j,\varphi_j)$ compactly contained in $D$ and containing $x_j$, and an open neighborhood $V_j$ of $x_j$ such that $V_j\Subset U_j$. Next, for each $j\in\{1,\ldots,m\}$,
 let $\chi_j\in\mathcal{C}^\infty_c(U_j)$ such that $\chi_j=1$ on $V_j$, and $\zeta_j(y):=|\varphi_j(y)|^2$ for $y\in U_j$. Then, set
 $\psi_j$ equal to $\chi_j\cdot\zeta_j$ on $U_j$ and zero on $M\setminus U_j$. Hence, each $\psi_j$, $j\in\{1,\ldots,m\}$, is a bounded, smooth function on $D$.
 Moreover, $i\partial\overline{\partial}\sum_{j=1}^m\psi_j\geq \omega$ on $\bigcup_{j=1}^mV_j$ for some volume form $\om$ on $M$. 
 Furthermore, $i\partial\overline{\partial}\psi_0\geq \omega$ 
 on $D\setminus\bigcup_{j=1}^mV_j$ for some volume form $\om$ on $M$.
 By the compactness of $\bigcup_{j=1}^m\overline U_j$ and the 
 smoothness of $\psi_j$, $j\in\{0,\ldots,m\}$, there exists an $\epsilon>0$ such that
\begin{align*}
  i\partial\overline{\partial}\left(\psi_0+\epsilon\sum_{j=1}^m\psi_j \right)\geq \omega\qquad{on}\;D
\end{align*}
for some volume form $\om$ on $M$. Therefore, $\psi:=\psi_0+\epsilon\sum_{j=1}^m\psi_j$ satisfies $(d)$ on $D$.

\subsection{Proof of $\mathbf{(d)}$ implies $\mathbf{(c)}$}\label{SS:dbar}
Let $(E,h)$ be a Hermitian holomorphic vector bundle over 
$M$, and set $(F,\mathfrak{h})=(T^{1,0}M\otimes E, 
g\otimes h)$. Then $(F,\mathfrak{h})$ is a Hermitian holomorphic vector bundle over $M$ such that
\begin{align*}
  A^2(D;E)\cong A^2(D,g;T^{*1,0}M\otimes T^{1,0}M
  \otimes E, g_{1,0}^*\otimes g\otimes h)
  \cong A^2_{1,0}(D;F).
\end{align*}
Hence, it suffices to show that $A_{1,0}^2(D;F)$ is infinite dimensional whenever ($d$) of Theorem \ref{T:main} holds. 

For fixed $N\in\mathbb{N}$, let $z_1,z_2,\ldots,z_N$ be distinct points in $D$. 
Choose charts $\left\{\left(U_j,\varphi_j\right)\right\}_{j=1}^N$, such that 
$U_j\subset D$, $z_j\in U_j$, and $\varphi_j(z_j)=0$. After 
possibly shrinking each $U_j$, we may assume that 
$U_j\cap U_k=\emptyset$ for any $j\neq k$ and 
$\varphi_j(U_j)$ is 
contained in the open unit disk $\mathbb{D}$ for each 
$j\in\{1,\ldots,N\}$.
Choose a function $\chi\in\mathcal{C}^\infty_c(\mathbb{D})$ 
which equals $1$ near the origin, and define 
$\ell_j\in\mathcal{C}^\infty_c(D)$, $j\in\{1,\ldots,N\}$, by 
$$\ell_j(z)=\chi(\varphi_j(z))\ln|\varphi_j(z)|\;\;\text{for}\;\;z\in U_j$$
and zero otherwise. Next, we set
$$
  \Phi_K:=K\psi+2\textstyle\sum_{j=1}^N\ell_j,
$$
where $\psi$ is the function supplied by the hypothesis, i.e.,
$\psi\in\mathcal{C}^\infty(D)$ is a bounded function which satisfies $i\partial\overline{\partial}\psi\geq\omega$ on $D$ for some volume form $\omega$ on $M$. 
It then follows that $\Phi_K\in\mathcal{C}^\infty(D_N)$ for $D_N:=D\setminus\{z_1,\ldots,z_N\}$ and any $K>0$. Furthermore, for any $\wt K>0$, there exists a 
$K>0$ such that
$\left(i\partial\dbar\Phi_{K}\right)\geq \wt K\omega_g$ on $D_N$. 
Note that $\mathfrak{h} e^{-\Phi_K}$ is a smooth Hermitian metric on $F$ over $D_N$, whose curvature is
\begin{align*}
 i\Theta(\mathfrak{h} e^{-\Phi_K})=i\Theta(\mathfrak{h})
 +i\partial\overline{\partial}\Phi_K\otimes \text{I}_F,
\end{align*}
on $D_N$, for any $K>0$.
Since $(F,\mathfrak{h})$ is a 
Hermitian vector bundle over $M$, there exists a constanct $c>0$ such that 
$i\Theta(\mathfrak{h})\geq -c\omega_g I_F$. We may now choose a $K>0$ such that $\bh:=\mathfrak h e^{-\Phi_K}$ satisfies  
$i\Theta(\bh)
\geq\omega_g I_{F}$ on $D_N$.
It follows that we may apply Theorem \ref{T:HDtheorem} for 
$(F,\bh)$ on $D_N$.

Next, let $\xi$ and $\eta$ be nonvanishing holomorphic sections of $F$ and $T^{1,0}M$, respectively, on $D$. Set
$\alpha=\overline{\partial}f\otimes\eta\otimes\xi$, where 
$f(z)=\chi(\varphi_N(z))$ for $z\in U_N$ and zero otherwise. It follows that $\alpha$ is a smooth, 
compactly supported $F$-valued 
$(1,1)$-form on $D_N$. Hence,
$\alpha\in L_{1,1}^2(D_N,g;F,\bh)$. Thus, by Theorem \ref{T:HDtheorem}, there exists a
$u\in L^2_{1,0}(D_N,g;F,\bh)$ such that $\dbar u=\alpha$ 
on $D_N$ in the sense of distributions and
\begin{align*}
  \langle u,u\rangle_{\bh}\leq
  \langle\alpha,\alpha\rangle_{\bh}.
\end{align*}
Then there exists a constant $c>0$ such that
\begin{align*}
  c\langle u,u\rangle_{\mathfrak{h}}
  \leq\langle \alpha,\alpha\rangle_{\bh}<\infty
\end{align*}
since $e^{-\Phi_K}$ has a positive lower bound on $D_N$.
Set $s=f\otimes\eta\otimes\xi-u$ on $D_N$. Then $\dbar s=0$ on $D_N$ 
in the sense of distributions. Thus,
$s\in\mathcal{O}(D_N; T^{1,0}M\otimes F)$. Moreover, since
both $f\otimes\eta\otimes\xi$ and $u$ belong to 
$ L^2_{1,0}(D_N,g;F,\mathfrak{h})$, so does $s$. It then follows that there exists an $\tilde s\in A^2_{1,0}(D,g;F,\mathfrak{h})$ such that $\tilde s|_{D_N}=s$, since $D\setminus D_N$ is a compact polar set; see \cite[Theorem 9.5]{Co95}. 
Therefore, $\tilde u:=f\otimes\eta\otimes\xi-\tilde s\in\Omega^{1,0}(D;F)$ is such that $\tilde{u}$ equals $u$ outside a set of measure zero, and $\langle \tilde{u},\tilde{u}\rangle_{\bh}=\langle u,u\rangle_{\bh}<\infty$. By the lack of local integrability of $e^{-\Phi_K}$ at $z_j$, we obtain that $\tilde{u}(z_j)=0$, $j\in\{1,\ldots,N\}$. Therefore,
$\tilde s(z_j)=0$ for all $j\in\{1,\ldots, N-1\}$ and $\tilde s(z_N)=\eta\otimes\xi$. As $N\in\mathbb{N}$ was arbitrary, it follows that $A^2_{1,0}(D,g;F,\mathfrak{h})$ is infinite dimensional.


\bibliography{referencesBergmanRS}{}
\bibliographystyle{plain}
\end{document}